\documentclass{article}
\newcommand{\R}{I\!\! R}
\newcommand{\N}{I\!\! N}
\newcommand{\Z}{I\!\! Z}
\newcommand{\C}{I\!\! C}

\begin{document}

\bigskip

\bigskip

Tsemo Aristide,

Visitor, Toronto University

100 Georges Street

Toronto

tsemoaristide@hotmail.com

\bigskip

\bigskip

\centerline{\bf An Atiyah-Singer Theorem for gerbes.}

\bigskip

\bigskip

\centerline{\bf Introduction.}

\bigskip

This paper has been motivated by the following problem: Let $M$ be
a compact riemannian manifold, the  curvatures are very useful to
study the topology of $M$. If the manifold $M$ is $spin$, the
study of the bundles of spinors provides results in this way.
Unfortunately every compact manifold is not spin. The obstruction
to the existence of a $spin-$structure on $M$ is the second
Stiefel-Whitney class $w_2(M)$ of $M$. Nevertheless the class
$w_2(M)$ is the classifying cocycle associated to a ${\Z}/2$ gerbe
on $M$ that we call the $spin$ gerbe, which is according to
Brylinski, and Mc Laughlin, an illuminating example of gerbe. The
objects of this gerbe are naturally endowed with a riemannian
metric, invariant by the automorphisms of the gerbe. It is natural
to think that the study of this $spin$ gerbe can have topological
applications. For example one may expect to generalize the
Lichnerowicz theorem. On this purpose, we need first to prove an
Atiyah-Singer type theorem for gerbes, which is our purpose.

\medskip

{\bf 1. On the notion of vectorial gerbes.}

\medskip

The aim of this section is to develop the notion of vectorial
gerbes.

\medskip

{\bf  Definition 1.}

 Let $M$ be a manifold, a sheaf $S$ of categories on $M$, is
a map $U\rightarrow S(U)$, where $U$ is an open set of $M$, and
$S(U)$ a category which satisfies the following properties:

- To each inclusion $U\rightarrow V$, there exists a map
$r_{U,V}:S(V)\rightarrow S(U)$ such that $r_{U,V}\circ
r_{V,W}=r_{U,W}$.

- Gluing conditions for objects,

Consider  a covering family $(U_i)_{i\in I}$ of an open set $U$ of
$B$, and for each $i$, an object $x_i$ of $S(U_i)$, suppose that
there exists a map $g_{ij}:r_{U_i\cap U_j,U_j}(x_j)\rightarrow
r_{U_i\cap U_j, U_i}(x_i)$ such that $g_{ij}g_{jk}=g_{ik}$, then
there exists an object $x$ of $C(U)$ such that $r_{U_i,U}(x)=x_i$

Gluing conditions for arrows,

Consider two objects $P$ and $Q$ of $S(M)$, then the map
$U\rightarrow Hom(r_{U,M}(P),r_{U.M}(Q))$ is a sheaf.

Moreover, if the following conditions are satisfied the sheaf of
categories $S$ is called a gerbe

$G1$

There exists a covering family $(U_i)_{i\in I}$ of $M$ such that
for each $i$ the category $S(U_i)$ is not empty

$G2$

Let $U$ be an open set of $M$, for each objects $x$ and $y$ of
$U$, there exists a covering family $(U_i)_{i\in I}$ of $U$ such
that $r_{U_i,U}(x)$ and $r_{U_i,U}(y)$ are isomorphic.

$G3$

Every arrow of $S(U)$ is invertible, and there exists a sheaf $A$
in groups on $M$, such that for each object $x$ of $S(U)$,
$Hom(x,x)= A(U)$, and the elements of this family of isomorphisms
commute with the restriction
 maps.

The sheaf $A$ is called the band of the gerbe $S$, in the sequel,
we will consider only gerbes with commutative band.

\medskip

{\bf Notation.}

\medskip

 For a covering family $(U_i)_{i\in I}$ of $B$, and an object $x_i$
 of $S(U_i)$, we denote
by $x^i_{i_1..i_n}$ the element $r_{({U_{i_1}}\cap..\cap
U_{i_n},U_i)}(x_i)$, and by $U_{i_1..i_n}$ the intersection
$U_{i_1}\cap..\cap U_{i_n}$.

\medskip

{\bf Definition 2.}

A gerbe is a vectorial gerbe if and only if for each open set $U$,
the category $S(U)$ is a category of vector bundles over $U$ with
typical fiber the vector space $V$ and maps between objects are
isomorphisms of vector bundles. The vector space $V$ will be
called the typical fiber of the vectorial gerbe. More precisely,
there exists a covering family $(U_i)_{i\in I}$ of $M$, a
commutative subgroup $H$ of $Gl(V)$, such that there exit maps
${g'}_{ij}:U_i\cap U_j\longrightarrow Gl(V)$, which define
isomorphisms
$$
g_{ij}: U_i\cap U_j\times V\longrightarrow U_i\cap U_j\times V
$$
$$
(x,y)\longrightarrow (x,g_{ij}(x)y)
$$
such that $c_{ijk}=g_{ij}g_{jk}g_{ki}$ is an $H$ $2-$Cech cocycle.

\medskip

{\bf Examples.}

\medskip

{\bf The Clifford gerbe associated to a riemannian structure.}

\medskip

Let $M$ be a $n-$riemannian manifold, $O(M)$ the reduction of the
bundle of linear frames which defines the riemannian structure of
the manifold $M$. The bundle $O(M)$ is a locally trivial principal
bundle over $M$ which typical fiber is $O(n)$. There is an exact
sequence $1\rightarrow {{\Z}/2}\rightarrow Spin(n)\rightarrow
O(n)\rightarrow 1$, where $Spin(n)$ is the universal cover of
$O(n)$. We can associate to this problem a gerbe which band is
${{\Z}/2}$ and such that for each open set $U$, $Spin(U)$ is the
category of $Spin(U)$ bundles over $U$ such that the quotient of
each of its element by ${{\Z}/2}$ is $O(U)$, the restriction of
$O(M)$ to $U$. The classifying cocycle of this gerbe is the second
Stiefel-Whitney class.

One can associate to this gerbe a vectorial gerbe named the
Clifford gerbe $Cl(M)$. For each open set $U$ of $M$, $Cl(U)$ is
the category which objects are Clifford bundles associated to the
objects of $Spin(U)$. The gerbe $Cl(M)$ is a vectorial gerbe.

Let ${g'}_{ij}$ be the transitions functions of the bundle $O(M)$,
for each map $g_{ij}$, consider an element ${g}_{ij}$ over
${g'}_{ij}$ in $Spin(n)$. The the element $g_{ij}(x)$ acts on
$Cl({\R}^n)$ by left multiplication, we will denote by $h_{ij}(x)$
the resulting automorphism of $Cl({\R}^n)$. The Clifford gerbe is
thus defined by $h_{ij}: U_i\cap U_j\rightarrow Spin(n)$.

\medskip

{\bf The gerbe defined by  the lifting problem associated to a
vectorial bundle.}

\medskip

Consider a vector bundle $E$ over $M$ which typical fiber is the
vector space $V$. One can associate to $E$, a principal $Gl(V)$
bundle. We suppose that this bundle has a reduction $E_K$ where
$K$ is a subgroup of $Gl(V)$. Consider a central extension
$1\rightarrow H\rightarrow G\rightarrow K\rightarrow 1$. This
central extension defines a gerbe $C_H$ on $M$, such that for each
open set $U$ of $M$, the objects of $C_H(U)$ are $G-$principal
bundles over $U$ which quotient by $H$ is the restriction of $E_K$
to $U$. We denote by $\pi$ the projection $\pi:G\rightarrow K$.

Suppose moreover defined a representation $r:G\rightarrow Gl(W)$,
and a surjection $f:W\rightarrow V$ such that the following square
is commutative
$$
\matrix{W &{r(h)\over\longrightarrow}& W\cr \downarrow f &\\& f
\downarrow\cr V&{\pi(h)\over \longrightarrow} &V}
$$

then one can defined the vectorial gerbe $C_{H,W}$ on $M$ such
that the object of $C(U)$ are $e_U\propto r$, where $e_U$ is an
object of $C_H(U)$. Let $(U_i)_{i\in I}$ be a trivialization of
$E$ defined by the transitions functions ${g'}_{ij}$, we consider
a map $g_{ij}:U_i\cap U_j\rightarrow G$ over ${g'}_{ij}$. The
gerbe $C_{H,W}$ is defined by $r(g_{ij})$.

Let $T$ be the set of elements  of $W$ fixed by elements of $H$.
The action of $G$ on $W$ defines an action of $K$ on $T$. This
action defines a vector bundle $S_T$ on $M$ with typical fiber is
$T$. Let $e_U$ be an object of the category $C_{H,W}$. then the
restriction of $S_T$ to $U$ is the set of elements of $e_U$
invariant by $H$.

\medskip

{\bf Definition 3.}

-A riemannian metric on a vectorial gerbe $C$ is defined by the
following data:

For each objects $e_U$ of $C$, a riemannian metric $<,>$ on the
vector bundle $e_U$ which is preserved by morphisms of between
objects of $C(U)$. We remark that the band need to be contained in
a compact group in this case since it preserves the riemannian
metric.

\medskip

An example of a scalar product on a gerbe is the following:
consider the Clifford gerbe $Cl(M)$, we know that the group $Spin$
is a compact group, its action on $Cl({\R}^n)$ preserves a scalar
product. This scalar product defines on each fiber of an object
$e_U$ of $Cl(U)$ a scalar product which defines, a riemannian
metric $<,>_{e_U}$ on $e_U$. The family of riemannian metrics
$<,>_{e_U}$ is a riemannian metric defined on the gerbe $Cl(M)$.

\medskip

{\bf Definition 4.}

 - A global section of a vectorial gerbe associated to a $1-$Cech
chain $(g_{ij})$ is defined by a covering space $(U_i)_{i\in I}$
of $M$, for each element $i$ of $I$, an object $e_i$ of $C(U_i)$,
a section $s_i$ of $e_i$, a family of morphisms
$g_{ij}:e^i_j\rightarrow e^j_i$ such that on $U_{ij}$ we have
$s_i=g_{ij}(s_j)$.

Let $(s_i)_{i\in I}$ be an element of $S(g_{ij})$, on $U_{ijk}$ we
have: $s_j=g_{jk}(s_k)$, $s_i=g_{ij}(s_j)$. This implies that
$s_i=g_{ij}g_{jk}(s_k)=g_{ik}(s_k)$. Or equivalently
$g^{-1}_{ik}g_{ij}g_{jk}(s_k)=s_k$. Remark that the restriction of
the element $s_k$ is not necessarily preserved by all the band.

Suppose that $M$ is compact, and $I$ is finite. We can suppose
that there exists $i_0$ such that $T=U_{i_0}-\cup_{i\neq i_0}U_i$
is not empty. Consider a section $s_{i_0}$ of $e_{i_0}$ which
support is contained in $T$, then we can define a global section
$(u_i)_{i\in I}$ such that $u_{i_0}=s_{i_0}$, and if $i\neq i_0$,
$u_i=0$. This ensures that $S(g_{ij})$ is not empty.

We will denote by $S(g_{ij})$ the family of global sections
associated to $(g_{ij})_{i,j\in I}$. Remark that $S(g_{ij})$ is a
vector space.

\bigskip

\medskip

{\bf Proposition 5.}

{\it Suppose that the vectorial gerbe $C$ is the gerbe associated
to the lifting problem defined by the extension $1\rightarrow
H\rightarrow G\rightarrow K\rightarrow 1$ and the vector bundle
$E$. Let $r:G\rightarrow Gl(W)$ be a representation, we suppose
that the condition  of the previous diagram  is satisfied. Then
each for $G-$chain $g_{ij}$, each element $(s_i)_{i\in I}$ of the
vector space of global sections $S(g_{ij})$, satisfies the
following condition: there exists a section $s$ of $E$, such that
$s_{\mid U_i}=f\circ s_i$.}

\medskip

{\bf Proof.}

Let $(s_i)_{i\in I}$ be a global section associated to the chain
$S(g_{ij})$, then on $U_{ij}$, we have $s_i=g_{ij}(s_j)$, this
implies that on $U_{ij}$, $f(s_i)=f(s_j)$. Thus the family
$(f(s_i))_{i\in I}$ of local sections of $E$ defines a global
section $s$ of $E$.

\medskip

{\bf Remark.}

\medskip

Let $s$ be a section of the bundle $E$, locally we can define a
family of sections $s_i$ of $e_i$, such that $f(s_i)=s_{\mid
U_i}$. We can consider the chain $s_{ij}=s_i-g_{ij}(s_j)$. We have
$s_{jk}-s_{ik}+s_{ij}=s_{ijk}$ is $2-$cocycle. Whenever there
exists a chain $g_{ij}$ a global section $s=(s_i)_{i\in I}$ such
that $s_i=g_{ij}(s_j)$, and $f(s_i)=s_{\mid U_i}$, it is not sure
that such a global section exists for another chain $h_{ij}$. This
motivates the following definition:

\medskip

{\bf Definition 6.}

We will define the vector space $S$ of formal global sections of
the vector gerbe $C$, the vector space which generators are $[s]$
where $s$ is an element of a set of global sections $S(g_{ij})$.

The elements of $S$, are formal finite sum of global sections.

\bigskip

{\bf The Prehilbertian structure of $S(g_{ij})$.}

\medskip

First we remark that $S(g_{ij})$ is a vector space. Let $s$, and
$t$ be elements of $S(g_{ij})$, we will denote by $s_i$ and $t_i$
the sections of $e_i$ which define respectively the global
sections $s$ and $t$. We have $s_i=g_{ij}(s_j)$ and
$t_i=g_{ij}(t_j)$ this implies that $as_i+bt_i=g_{ij}(as_j+bt_j)$
where $a$ and $b$ are real numbers.

\medskip

- The scalar structure of $S(g_{ij})$.

\medskip

Let $(V_k, f_k)_{k\in K}$ be a partition of unity subordinate to
$(U_i)_{i\in I}$, this means that for each $k$ there exists an
$i(k)$ such that $V_k$ is a subset of $U_{i(k)}$. Since the
support of $f_k$ is a compact subset of $V_k$, we can calculate
$\int<s_{ii(k)},t_{ii(k)}>$ where $s_{ii(k)}$ and $t_{ii(k)}$ are
the respective restrictions of $s_{i(k)}$ and $t_{i(k)}$ to $V_k$,
remark that since we have supposed that $s_i=g_{ij}(s_j)$ and
$g_{ij}$ is an riemannian isomorphism between $e_j^i$ and $e_i^j$,
if $V_k$ is also included in $U_j$, then $<s_{ii(k)},t_{ii(k)}>=
<s_{jj(k)},t_{jj(k)}>$ on $V_k$. We can define
$<s,t>=\sum_k\int<f_ks_{ii(k)},f_kt_{ii(k)}>$.

We will denote by $L ^2(S(g_{ij})$ the Hilbert completion of the
Pre Hilbert structure of $(S(g_{ij}),<,>)$.

\medskip

{\bf The scalar structure on set of formal global sections $S$.}

\medskip

Let $s$ and $t$ be two formal global sections, we have
$s=[s_{n_1}]+..+[s_¸{n_p}]$, and $t=[t_{m_1}]+..+[t_{m_q}]$, where
$s_{n_i}$ and $t_{m_j}$ are global sections.

We will define a scalar product on $S$ as follows: if $s$ and $t$
are elements of the same set of global sections $S(g_{ij})$,
$<[s],[t]>=<s,t>_{S(g_{ij})}$. If $s$ and $t$ are not elements of
the same set of global sections, then $<[s],[t]>=0$.

\medskip

{\bf Proposition 7.}

{\it An element of $L^2(S(g_{ij}))$ is a family of $L^2$ sections
$s_i$ of $e_i$ such that $s_i=g_{ij}(s_j)$.}

\medskip

{\bf Proof.}

Let $(s^l)_{l\in N}$ be a Cauchy sequence of $(S(g_{ij}),<,>)$. We
can suppose that the open sets $V_k$ used to construct the
riemannian metric are such that the restriction $e^k_i$ of $e_i$
to $V_k$ is a trivial vector bundle. The sequence
$(f_ks^l_{ii(k)})$ is a Cauchy sequence defined on on the support
$T_k$ of $f_k$. Since this support is compact, we obtain that
$(f_ks^l_{ii(k)})_{l\in N}$ goes to an $L^2$ section $s_{i(k)}$ of
$e^k_i$. We can define $s_i=\sum_{k,V_k\cap
U_i\neq\phi}f_ks_{i(k)}$. The family $(s_i)_{i\in I}$ defines the
requested limit.

\medskip

Suppose that morphisms between objects commute with laplacian, we
can then endow $S(g_{ij})$ with the prehilbertian structure
defined by $<u,v>=\int<\Delta^s(u),v>$, where $\Delta^s(u)$ is the
global section defined by ${\Delta^s(u)}_i=\Delta^s(u_i)$. We will
denote by $H_s(S(g_{ij}))$ the Hilbert completion of this
prehilbertian space.

We will define the formal $s-$distributional global sections
$H_s(S)$ as the vector space generated by finite sums
$[s_1]+..[s_k]$ where $s_i$ is an element of an Hilbert space
$H_s(S(g_{ij}))$.

\medskip

{\bf Connection on riemannian gerbes and characteristic classes.}

\medskip

The notion of connection is not well-defined for general vectorial
gerbe, nevertheless the existence of a riemannian structure on a
vectorial gerbe, $C$, gives rise to a riemannian connection on
each object $e_U$ of $C(U)$, this family of riemannian connections
will be the riemannian connection of the gerbe $C$.

Let $(U_i)_{i\in I}$ be an open covering of $M$, suppose that the
objects of $C(U_i)$ are trivial bundles. The riemannian connection
of the object $e_i$ of $C(e_i)$ is defined by a $1-$ $so(V)$ form
$w_i$ on $TU_i$, and the covariant derivative of this connection
evaluated to a section $s_i$ of $e_i$ is $ds_i+w_is_i$. The
curvature of this connection is the $2-$form
$\Omega_i=dw_i+w_i\wedge w_i$.

The $2k^{eme}$ Chern class of $e_i$, ${c^i}_{2k}$ is defined by
$Trace[({i\over 2\pi}\Omega_i)^k]$. Let $e'_i$ be another object
of $C(U_i)$. There exist isomorphisms $\phi_i:e_i\rightarrow
U_i\times V$, and $\phi_i':e'_i\rightarrow U_i\times V$., and
$g_i:e_i\rightarrow e'_i$. The map $\phi_i'\circ g_i\circ
\phi_i^{-1}$ is an automorphism of $U_i\times V$ defined as
follows:
$$
 h_i:(x,y)\rightarrow (x,u_i(x)y))
$$
The riemannian ${\phi_i}^{-1}<,>={\phi'_i}^{-1}<,>$ is preserved
by $h_i$. We have
${{\phi_i}^{-1}}^*\Omega_i=u_i{{\phi_i'}^{-1}}^*{\Omega'_i}{{u_i}^{-1}}$.
This implies that ${c^i}_{2k}={{c'}^i}_{2k}$.

\medskip

 There exists an isomorphism between the respective restrictions
 $e^i_j$ and $e^j_i$ on $U_i$. As above we can show that this implies
 that the $2k-$Chern classes of $e^i_j$ and $e^j_i$ coincide  on $U_i\cap
 U_j$, and defines a global class $c_{2k}$ on $M$ which is the $2k-$ Chern
 class of the riemannian gerbe.

 We can define the $c(C)=c_1(M)+...+c_n(M)$ the total Chern form
 of the gerbe, and the form define locally by ${ch(C)}_{\mid U_i}=
 Tr(exp(i{\Omega\over 2\pi}))$ the total Chern character.

\medskip

{\bf 2. Operators on riemannian gerbe.}

\medskip

We begin by recalling the definition of pseudo-differential
operators for open sets of ${\R}^n$ and manifolds.

 Let $U$, be an open set of
${\R}^n$, we denote by $S^m(U)$ the set of smooth functions
$p(x,u)$ defined on $U\times{\R}^n$, such that for every compact
set $K\subset U$, and every multi-indices $\alpha$ and $\beta$, we
have $\mid\mid
D^{\alpha}D^{\beta}p(x,u)\mid\mid<C_{\alpha,\beta,K}(1+\mid\mid
u\mid\mid)^{m-\mid \alpha\mid}$.

Let $K(U)$ and $L(U)$ denote respectively the smooth functions
with compact support defined on $U$ and the smooth functions on
$U$. We can define the map, $P:K(U)\rightarrow L(U)$ such that
$$
P(f)=\int p(x,u)\hat f(u)e^{i<x,u>}du
$$
where $\hat f$ is the Fourier transform of $f$.

\medskip

{\bf Definition 1.}

 An operator on $U$  is pseudo-differential, if it is locally of
the above type.

\medskip

{\bf Definition 2.}

Let $P$ be a pseudo-differential operator, $(U_i)_{i\in I}$ a
covering family of $U$ such that the restriction of $P$ to $U_i$
is defined by $P(f)=\int p_i(x,u)\hat f(u)e^{i<x,u>}du$. The
operator is of degree $m$ if
$\sigma(p_{\mid{U_i}})=lim_{t\rightarrow\infty}{p_{i}(x,tu)\over
t^m}$ exists. In this case $\sigma(p)$ is called the symbol of
$P$.

\medskip

Let $E$ be a vector bundle over the riemannian manifold $M$,
endowed with a scalar metric. We denote by $K(E)$ and $L(E)$ the
respectively  set of smooth sections of $E$ with compact support
and the set of smooth sections of $E$. An operator on the vector
bundle $E$, is a map $P:K(E)\rightarrow L(E)$ such that there
exists a covering family $(U_i)_{i\in I}$ which satisfies:

- The restriction of $E$ to $U_i$ is trivial

- The restriction of $P_i$ of $P$ to $U_i$ is a map
$P_i:K(U_i\times V)\rightarrow L(U'_i\times V)$ where $V$ is the
typical fiber of $E$.

- If we consider charts $\phi_i$ and $\psi_i$ such that
$\phi_i(U_i\times V)=\psi_i(U'_i\times V)=U\times {\R}^n$, then
the map $P_i$ is defined by the a matrix $(p_{kl})$ where $p_{kl}$
 define an operator of degree $m$. More precisely, if $s'$ is a section
  of $E$ over $U_i$
and $s=(s_1,..,s_n)=\phi_i(s')$, we can define
$t_k=\sum_{l=1}^{l=n}\int p_{lk}(x,u)\hat s_l(u)e^{i<x,u>}du$, and
$P_i(s')={\psi_i}^{-1}(t_1,..,t_n)$.

 Consider $SM$ the sphere bundle of the cotangent space $T^*M$ of $M$,
and $\pi^*E$ the pull-back of $E$ to $T^*M$, the symbols defined
by $(p_{ij})$ define a map $\sigma:\pi^*E\rightarrow \pi^*E$.
Consider now the projection $\pi_S:SM\rightarrow M$, then $\sigma$
induces a map $\sigma_S:{\pi_S}^*E\rightarrow {\pi_S}^*E$.

 Let $s$ be a positive integer
 we denote by ${{H_s}^{loc}}(M,E)$, the space of distributions sections $u$ of $E$
such that $D(u)$ is a ${L^2}^{loc}$ section, where $D$ is any
differential operator of order less than $s$, and by
${{H_s}^{comp}}(M,E)$ the subset of elements of
${{H_s}^{loc}}(M,E)$ with compact support. Remark that if $M$ is
compact, then ${{H_s}^{loc}}(M,E)={{H_s}^{comp}}(M,E)$. We define
by ${H_{-s}}^{loc}(M,E)$ to be the dual space of
${H_s}^{comp}(M,E)$, and by ${H_{-s}}^{comp}(M,E)$ the dual space
of ${H_s}^{loc}(M,E)$.

The Sobolev space $H_s$ is an Hilbert space endowed with the norm
defined by $(\mid\mid \int <\Delta^su,u>\mid\mid)^{1\over 2}$.

Every operator $P$ of order less than $m$ can be extended to a
continuous morphism $H_s\rightarrow H_{s-m}$.

\medskip

{\bf Definition 3.}

Let $C$ be a riemannian gerbe defined on the manifold $M$, an
operator $D$ of degree $m$ on $C$, is a family of operators
$D_{e}$ of degree$m$ defined on $e$ where $e$ is an object of the
category $C(U)$. We suppose that for each morphism $g:e\rightarrow
f$, $D_fg^*=g^*D_e$.

\medskip

{\bf Remark.}

The last condition in the previous definition  implies that $D_e$
is invariant by the automorphisms of $e$. The map $g^*$  is the
map which transforms a section $s$ to $g(s)$. The operator
considered in the sequel will be assumed to be continue. It
defines a map $D_e:{{H_s}^{comp}}(U,e)\rightarrow
{H_{s-m}^{loc}}(U,e)$.

\medskip

{\bf Proposition 4.} {\it Let $D$ be an operator of degree $m$
defined on the riemannian gerbe $C$, then $D$ induces a map
$D_{S(g_{ij})}:H_s(S(g_{ij}))\rightarrow H_{s-m}(S(g_{ij}))$ and a
map $D_S:H_s(S)\rightarrow H_{s-m}(S)$.}

\medskip

{\bf Proof.}

Consider a global distributional section $s$ which is an element
of $H_s(S_{g_{ij}}))$, we have
$g_{ij}(D_{e^i_j}(s^i_j))=D_{e^j_i}(s^j_i)$. This implies the
result.

\medskip

In the sequel we will consider only pseudo-differential operators
that preserve $C^{\infty}$ sections.

\bigskip

{\bf The symbol of an operator.}

\bigskip

Let $C$ be a vectorial gerbe defined on $M$ endowed with the
operator $D$ of degree $m$, for each object $e$ of $C(U)$, we can
pulls back the bundle $e$ by the projection map
$\pi_{SU}:SU\rightarrow U$ to a bundle $\pi_{SU}^*e$ over $SU$,
where $SU$ is the restriction of the cosphere bundle defined by a
fixed riemannian metric of $T^*M$. The family $C_S(U)$ which
elements are $\pi_{SU}^*e$ is a category. The maps between objects
of this category are induced by maps between elements of $C(U)$.
The map $U\rightarrow C_S(U)$ is a gerbe which has the same band
than $C$.

Now on the object $e$, we can define symbol
$\sigma_{D_e}:\pi_{SU}^*e\rightarrow \pi_{SU}^*e$.

Remark that for every automorphism $g$ of $e$, the fact that
$g^*\circ D_e=D_e\circ g^*$ implies that
$\sigma_{gD_eg^{-1}}=\sigma_{D_e}$.

\bigskip

\medskip

{\bf Proposition 5. Rellich Lemma for the family $S(g_{ij})$.}

{\it Let $f_n$ be a sequence of elements of $H_s(S(g_{ij}))$, we
suppose that there is a constant $L$ such that $\mid\mid
f_n\mid\mid_s<L$, then for every $s>t$,  there exits a subsequence
$f_{n_k}$ which converges in $H_t$.}

\medskip

{\bf Proof.}

Let $(s_n)$ be a sequence of sections which satisfy the condition
of the proposition, and $(V_{\alpha}, f_{\alpha})$ a partition of
unity subordinate to $(U_i)_{i\in I}$. We suppose that the support
of $f_{\alpha}$ is a compact space $K_{\alpha}$. We denote by
$s^i_n$ the section of $e_i$ which defines $s_n$, and by
$s^i_{n\alpha}$ the restriction of $s^i_n$ to the restriction of
$e_i$ to $V_{\alpha}$. The family $(f_{\alpha}s^i_{n\alpha})$ goes
to the element $s_{i\alpha}$ in $V_{\alpha}$ by the classical
Rellich lemma. We can write then $s^i=\sum_{V_{\alpha}\subset
U_i}f_{\alpha}s_{i\alpha}$. this is an $H_t$ map since the family
of $V_{\alpha}$ can be supposed to be finite, since $M$ is
compact. The family $(s^i)$ define a global $H_s$ section which is
the requested limit.

\medskip

{\bf Remark.}

A compact operator between Hilbert space is an operator which
transforms bounded spaces to compact spaces. The previous lemma
implies that if $s>t$, then the inclusion
$H_s(S(g_{ij}))\rightarrow H_t(S(g_{ij}))$ is compact, since as
$H_s(S(g_{ij}))$ is a separate space, a compact subspace of
$H_s(S(g_{ij}))$ is a set such that we can extract a convergent
sequence from every bounded sequence.

\medskip

{\bf Proposition 6.}

{\it The space $Op(C)$  of continuous linear maps of
$H_s(S(g_{ij}))$ is a Banach space.}

\medskip

{\bf Proof.}

Let $D_n$ be a Cauchy sequence of elements of $Op(C)$, for each
global section $s$, the sequence $D_n(s)$  is a Cauchy sequence in
respect to the norm of $H_s$, we conclude that it goes to an
element $D(s)$. The map $D:s\rightarrow D(s)$ is the requested
limit. It is bounded since $(D_n)_{n\in N}$ is a Cauchy sequence.

\medskip

The previous Lemma allows us to define $O^m$, the completion of
the pseudo-differential operators in   $OP(C)$ of order $m$, and
to extend the symbol $\sigma$ to $O^m$.

Now we will show that  the kernel of the extension of the symbol
to $O^m$ contains only compact operators.

\medskip

{\bf Definition 7.}

We say that an operator is elliptic if the family of symbol
$\sigma_{D_e}$ are invertible maps.

\medskip

{\bf Proposition 8.}

{\it The kernel of the f of the symbol map contains only compact
operators.}

\medskip

{\bf Proof.}

The symbol $\sigma(P)$ of the operator $P$ is zero if and only if
the order $m$ of the operator is  less than $-1$. This implies
that the operator $P$ is compact. To see it, we can suppose our
operator to be an $L^2(S(g_{ij}))$ operator, by composing it by
the inclusion map $H_{2-s}(S(g_{ij}))\rightarrow L^2(S(g_{ij})$,
we conclude by using the previous Rellich lemma.

\medskip

{\bf 3. $K-$theory and the index.}

\medskip

In this part we will give the definitions of the $K-$theory groups
$K_0$, and $K_1$, and show how we can use them to associate to a
symbol of an operator on a gerbe, an element of $K_0(T^*M)$.

\medskip

We will denote by $M_n$ the vector space of $n\times n$ complex
matrix. For $n\leq m$, we consider the natural injection
$M_n\rightarrow M_m$. We will call $M_{\infty}$, the inductive
limit of the vector spaces $M_n, n\in{\N}$.

Let $R$ be a ring and $p$, and $q$ be two idempotents of
$R_{\infty}=R\otimes M_{\infty}$, we will say that $p\sim q$ if
and only if there exists elements $u$ and $v$ of $R_{\infty}$ such
that $p=uv$, and $q=vu$. We denote by $[p]$ the class of $p$, and
by $Idem(R_{\infty})$ the set of equivalence classes.

 If $[p]$ and $[q]$ are represented respectively by elements of $R\otimes M_n$ and
$R\otimes M_m$, we can define an idempotent of $R\otimes M_{n+m}$
represented by the matrix $\pmatrix{p&0\cr 0&q}=[p+q].$

\medskip

{\bf Definition 1.}

We will denote by $K_0(R)$, the semi-group $Idem(R_{\infty})$,
endowed with the law $[p]+[q]=[p+q]$.

\medskip

Let $M$ be a compact manifold, and $C(M)$ the set of complex
valued functions on $M$. It is a well-known fact that for a
complex vector bundle $V$ on $M$, there exists  a bundle $W$, such
that $V\oplus W$ is a trivial bundle  isomorphic to $M\times
{\C}^k$. We can thus identify a vector bundle over $M$ to an
idempotent of $C(X)\otimes M_k$ which is also an idempotent of
$C(X)_{\infty}$. This enables to identify $K_0(M)$ to $K_0(C(X))$.

In fact the semi-group $K(R)$ is a group.

\medskip

Let $Gl_k(R)$ be the group of invertible elements of $M_k(R)$, if
$l\leq k$ we have the canonical inclusion map $Gl_l(R)\rightarrow
Gl_k(R)$. We will denote by $Gl_{\infty}(R)$, the inductive limit
of the groups $Gl_k(R)$.

\medskip

{\bf Definition 2.}

Let $Gl_{\infty}(R)_0$ be the connected component of
$Gl_{\infty}(R)$. We will denote by $K_1(R)$ the quotient of
$Gl_{\infty}(R)$ by $Gl_{\infty}(R)_0$.

\medskip

For a compact manifold $M$, we define $K_1(M)$ by $K_1(C(X))$.

Consider now an exact sequence $0\rightarrow R_1\rightarrow
R_2\rightarrow R_3\rightarrow 0$ of $C^*$ algebras, we have the
following exact sequences in $K-$theory:
$$
K_1(R_1)\rightarrow K_1(R_2)\rightarrow K_1(R_3)\rightarrow
K_0(R_1)\rightarrow K_0(R_2)\rightarrow K_0(R_3)
$$

\medskip

Let ${\cal H}$ be an Hilbert space, we denote by $B({\cal H})$ the
space of continuous operators defined on ${\cal H}$, and ${\cal
K}$ the subspace of compact continuous operators. We have the
following exact sequence:
$$
0\rightarrow {\cal K}\rightarrow B({\cal H})\rightarrow B({\cal
H})/{\cal K}=Ca\rightarrow 0.
$$
It is a well-known fact that $K_0({\cal K})={\Z}$.

\medskip

 Let $M$ be a riemannian manifold, and $C$ a riemannian gerbe
defined on $M$. Consider an elliptic operator $D$ defined on $C$,
of degree $l$. The operator $D$ induces a morphism:
$D:L^2(S(g_{ij}))\rightarrow H_{2-l}(S(g_{ij}))$. Consider the
operator  $(1-\Delta)^{-m}$ of degree $-l$. The operator
$D'=(1-\Delta)^{-m}D$ is a morphism of $L^2(S(g_{ij}))$.  The
symbol of $(1-\Delta)^{-m}D=\sigma(D)$. This implies that the
image of the operator $D$ in the Calkin algebra $L^2(S(g_{ij}))$
is invertible. It thus define a class $[\sigma(D')]$ of $K_1(Ca)$.
The image of $[\sigma(D')]$ in
 $K_0({\cal K})$ is the index of $D$. We remark that
the index of the operator depends only of the symbol.

For every object $e$ of $C(U)$, the symbol $\sigma({D_e})$ is an
automorphism of $\pi_{SU}e$, it defines an element $[\sigma(D_e)]$
of $K_1(C(S(U)))$ (recall that $S(U)$ is the cosphere bundle over
$U$ defined by the riemannian metric).

\bigskip

{\bf Remark.}

Let $U$ be an open set such that the objects of $C(U)$ are trivial
bundles. Consider an object $e$ of $C(U)$, and a trivialization
map $\phi_e:e\rightarrow U\times V$. For every object $f$ of
$C(U)$, we have
${{\phi_e}^{-1}}^*(\sigma_{D_e})={{\phi_f}^{-1}}^*(\sigma_{D_f})$.

\medskip

{\bf Proposition 3.}

{\it Let $C$ be a vectorial gerbe defined over the compact
manifold $M$, there is a trivial complex bundle $f_m=M\times
{\C}^m$ such that each object of $C(U)$ where $U$ is an open set
of $M$ is isomorphic to  a subbundle of the restriction of $f_m$
to $U$.}

\medskip

{\bf Proof.}

Let $(U_i)_{i\in I}$ be a finite covering family such that for
each $i$ the objects of the category $C(U_i)$ are trivial bundles.
Let $e_U$ be an object of $C(U)$, the restriction $e_i$ of $e$ to
$U_i$ is a trivial vectorial bundle. We consider a fixed  object
$e_i^0$ of $C(U_i)$. Consider a finite partition of unity
$(f_p)_{p\in 1,..,l}$ subordinate to the covering family $(U_i)$,
and $g_i:e_i^0\rightarrow V$ the composition of the trivialization
and the second projection. Let $h_i:e_i\rightarrow e_i^0$, we can
define the map $k:e_U\rightarrow {\C}^{ldimV} $ such that
$k(x)=(f_1(\pi_{e_U}(x))g_1h_1(x),..,f_l(\pi_{e_U}(x))g_lh_l(x))$.
$k$ induces a map $K:M\rightarrow G_m({\C}^N)$, where $m=dim V$
and $N=ldim(V)$. $e_U$ is the pulls-back of the canonical
$m-$vector bundle over $G_{m}({\C}^N)$. We remark that its a
subbundle of the pulls-back of the trivial bundle
$G_m({\C}^N)\times {\C}^N$.

\medskip

The orthogonal bundle of $e_U$ in the previous result can be
chosen canonically by considering the orthogonal bundle of the
canonical ${\C}^m$ bundle over $G_m({\C}^N)$ in $G_m({\C}^N)\times
{\C}^N$. We will suppose that this bundle is chosen canonically in
the sequel.

\medskip

 {\bf Proposition 4.}

{\it Let $C$ be a riemannian gerbe defined over the compact
manifold $M$. Then we can associate naturally to the symbol of the
elliptic operator $D$ a class $[\sigma_D]$ in $K_1(T^*M)$.}

\medskip

{\bf Proof.}

Let $(U_i)_{i\in I}$, a finite covering family of $M$ such that
for each $i$, each object of $e_i$ is a trivial bundle. Consider
for each $i$ a trivialization $\phi_i:e_i\rightarrow U_i\times V$,
The map ${\phi^{-1}}^*(\sigma_{D_{e_i}})$ is an automorphism of
the bundle ${\phi_i^{-1}}^*({\pi_S}^*(e))=S_e(U_i)$. This enables
to extend ${\phi^{-1}}^*(\sigma_{D_{e_i}})$ to a morphism of the
restriction of $f_m$ to $U_i$ completing by $1$ on the diagonal
since the orthogonal of $e_U$ is chosen canonically. It results a
morphism $\sigma_D'$ of $C(M)\otimes {\C}^m$, that is as an
element of $C(S^*M)\otimes Gl(V)$.

This element defines the requested element $[\sigma_D']$ of
$K_1(C(S^*M)\otimes M_k)\simeq K_1(C(S^*M))$.

\medskip

Let $B^*M$ be the compactification of $T^*M$ with fibers
isomorphic to the unit ball. We consider the bundle $B^*M/T^*M$
that we can identify to the sphere bundle $S^*M$. We have the
following exact sequence
$$0\rightarrow C_0(T^*M)\rightarrow
C(B^*M)\rightarrow C(S^*M)\rightarrow 0.
$$
This sequence gives rise to the following exact sequence in
$K-$theory:
$$
K_1(C(S^*M)\otimes M_k)\rightarrow K_0(C(T^*M)\otimes
M_k)\rightarrow K_0(C(B^*M)\otimes M_k)\rightarrow
K_0(C(S^*M)\otimes M_k)\rightarrow 0.
$$
We can define the boundary $\delta([\sigma_D'])$ which is an
element of $K_0(T^*M\otimes M_k)\simeq K_0(T^*M)$.

\medskip

{\bf Proposition 5.}

{\it The index of $D$ depends only of the class of
$\delta[\sigma_D']$ in $K_0(T^*M)$.}

\medskip

{\bf Proof.}

We remark that the kernel of the image of a symbol by the map
$K_1(S^*M\otimes M_k)\rightarrow K_0(T^*M\otimes M_k)$ is zero, if
it is the restriction of a map defines on $B^*M$. This implies
that this symbol is homotopic to a function which does not depend
of the second variable $u$. Thus the operator it defines is
homotopic to a multiplication by a constant which index is zero.

\bigskip

E. Meinrenken asked me the following question:

Let $M$ be a manifold which is the union of two open sets $U_1$
and $U_2$ such that there exists objects $e_1$ and $e_2$ of the
respective categories $C(U_1)$ and $C(U_2)$. Given operators
$D_{e_1}$ and $D_{e_2}$ on $e_1$ and $e_2$ is it possible to
associate to these operators an element of the $K-$theory? We do
not request any compatibility between $D_{e_1}$ and $D_{e_2}$.

\medskip

 The
vectors bundles $e_1$ and $e_2$ are subbundles of the respective
trivial bundles $U_1\times {\C}^k$ and $U_2\times {\C}^k$. Let
$SU_1$ and $SU_2$ be the restriction of the sphere bundle of the
cotangent space of $M$ to $U_1$ and $U_2$, we can canonically
extends the symbols of $D_{e_1}$ and $D_{e_2}$, to respective
automorphisms $\sigma_{D_1}$ and $\sigma_{D_2}$ of $F_1=SU_1\times
{\C}^k$, and $F_2=SU_2\times {\C}^k$, thus elements of
$K_1(C(U_1)\otimes M_k)=K_1(C(U_1))$, and $K_1(C(U_2)\otimes
M_k)=K_1(C(U_2))$, where $C(U_1)$ and $C(U_2)$ are respectively
the set of differentiable functions of $U_1$ and $U_2$.

 The Mayer-Vietoris sequence for algebraic $K-$theory gives rise
 to the sequence
$$
K_2(C(S(U_1\cap U_2))\rightarrow K_1(C(SM))\rightarrow
K_1(C(SU_1))\oplus K_1(C(SU_2))\rightarrow K_1(C((S(U_1\cap U_2)))
$$
if the image of $[\sigma_{D_{e_1}}]+[\sigma_{D_{e_2}}]$ in the
previous sequence is zero, then there exists a class $[\sigma_D]$
of $K_1(SM)$, which image by the map of the previous exact
sequence is the element $[\sigma_{D_{e_1}}]+[\sigma_{D_{e_2}}]$ of
$K_1(C(SU_1)\oplus K_1(C(SU_2))$. The class $[\sigma_D]$ is not
necessarily unique.

\medskip

\bigskip

{\bf The index formula for operator on gerbes.}

\medskip

Now, we will deduce an Atiyah-Singer type theorem for riemannian
gerbe. We know that the Chern character of the cotangent bundle
induces an isomorphism:
$$
K_0(T^*M)\otimes{\R}\rightarrow H_c^{even}(M,{\R})
$$
$$
x\otimes t\rightarrow tch(x)
$$
Consider $Vect(Ind)$ the subspace of $K_0(T^*M)$ generated by
$\sigma_P$, where $P$ is an operator on the riemannian gerbe. It
can be considered as a subspace of $H_c^{even}(M,{\R})$. The map
$Vect(\sigma_P)\rightarrow {\R}$ determined by $ch([\sigma
P'])\rightarrow ind(P)$ can be extended to a linear map
${H_c}^{even}(M,{\R})\rightarrow {\R}$.

The Poincare duality implies the existence of a class $t(M)$ such
that
$$
Ind(P)=\int_{T^*M} ch([\sigma_P'])\wedge t(M)
$$

\bigskip

{\bf 4. Applications.}

\bigskip

We will apply now this theory to the problem that has motivated is
construction.

Let $M$ be a riemannian manifold, consider the  Clifford gerbe on
$M$,

Let $(U_i)_{i\in I}$ be an open covering of $M$,  the riemannian
connection $w_{Civ}$ is defined by a family of $so(n)$ forms $w_i$
on $U_i$ which satisfy
$w_j={adg_{ij}}^{-1}w_i+{g_{ij}}^{-1}dg_{ij}$. The covariant
derivative of the Levi-Civita connection is $d+w_i$. We will fixes
an orthogonal basis $(e_1,..,e_n)$ of the tangent space $TU_i$ of
$U_i$. and write $w_i=\sum_{k=1}^{k=n}w_{ik}e_k$.

 we can define the spinorial covariant derivative by setting
$\phi_{ij}=-{1\over 4}w_{ij}$. In the orthogonal basis
$(e_1,..,e_n)$ we have
$\phi_{e_j}=\sum_{k,l}{\phi_{kl}}_{kl}e_ke_l$, with
${\phi}_{kl}=-{\phi}_{lk}$.

\medskip

{\bf The Dirac operator.}

\medskip

Let $e_U$ be an object of $Cl(U)$, $U_i$ a trivialization of
$Cl(U)$, we will define
$D_{e_U}=\sum_{k=1}^{k=n}e_i{\nabla_{spin}}_{e_i}$,

On each object $e_U$ of $C(U)$, we have the
Lichnerowicz-Weitzenbock formula: $D^2=\nabla^*\nabla+{1\over
4}s$, where $s$ is the scalar curvature. In this formula
$\nabla^*\nabla$ is the connection laplacian.

We will say the global spinor is harmonic if $D_{e_i}(s_i)=0$, for
each $s_i$.

\medskip

{\bf Proposition 1.}

{\it Suppose that the scalar curvature $s$ of $M$  is strictly
positive, and $M$ is compact then every harmonic global spinor is
$0$.}

\medskip

{\bf Proof.}

Let $\psi$ be an harmonic global spinor, we can represent $\psi$
by a family of spinor $s_i$ defined on an open cover $(U_i)_{i\in
I}$ of $M$, we have on each $U_i$, $D_{e_i}(s_i)=0$, this implies
that ${D_{e_i}}^2(s_i)=0$, we can write
$$
\int_{U_i}\nabla\nabla^*s_i+{1\over 4}s=0
$$
this implies that $\int_{M}s=0$, which  contradicts the fact that
the scalar curvature is strictly positive.

\medskip

\medskip

{\bf Corollary 2.} {\it Suppose that  the sectional curvature of a
compact riemannian manifold is strictly  positive  then the class
$\tau(M)$ associated to the index formula for operators on the
$Cl(M)$ gerbe is zero.}

\bigskip

{\bf Acknowledgements.} The author would like to thanks E.
Meinrenken for helpful comments.

\bigskip

{\bf Bibliography.}

\bigskip

1. Atiyah, M. Singer, I. The index of elliptic operators, I. Ann.
of Math (87) 1968, 484-530.

\medskip

2. Brylinski, J. Mc Laughin The geometry of degree-four
characteristic classes and line bundles on loop spaces I. Duke
Math. Journal (75) 603-637

\medskip

3. Charles, J. Cours de maitrise de mathematique

\medskip

 4. Van Erp, E. The Atiyah-Singer index theorem,
 $C^*$ algebraic $K-$theory and  quantization
 Master thesis, University of Amsterdam

\medskip

5. Tsemo, A. The differentiable geometry of discrete vectorial
gerbes, in preparation.

\end{document}